# PRACTICAL DRIFT CONDITIONS FOR SUBGEOMETRIC RATES OF CONVERGENCE

By Randal Douc, Gersende Fort, Eric Moulines and Philippe Soulier

*Ecole Polytechnique, LMC-IMAG, Ecole Nationale Supérieure des Télécommunications and Université d, Evry*

We present a new drift condition which implies rates of convergence to the stationary distribution of the iterates of a $\psi$-irreducible aperiodic and positive recurrent transition kernel. This condition, extending a condition introduced by Jarner and Roberts [*Ann. Appl. Probab.* **12** (2002) 224–247] for polynomial convergence rates, turns out to be very convenient to prove subgeometric rates of convergence. Several applications are presented including nonlinear autoregressive models, stochastic unit root models and multidimensional random walk Hastings–Metropolis algorithms.

**1. Introduction.** Let $(\Phi_n, n \geq 0)$ be a discrete time Markov chain on a general measurable state space $(\mathsf{X}, \mathcal{B}(\mathsf{X}))$ with transition kernel $P$. Assume that it is $\psi$-irreducible, aperiodic and positive recurrent. This paper considers the use of drift conditions to establish the convergence in $f$-norm of the iterates $P^n$ of the kernel to the stationary distribution $\pi$ at rate $r := (r(n), n \geq 0)$; that is,

$$(1.1) \qquad \lim_{n \to \infty} r(n) \| P^n(x, \cdot) - \pi \|_f = 0, \qquad \pi \text{ a.e.}$$

where $f : \mathsf{X} \to [1, \infty)$ satisfies $\pi(f) < \infty$ and for any signed measure $\mu$, the $f$-norm $\|\mu\|_f$ is defined as $\sup_{|g| \leq f} |\mu(g)|$.

For geometric rate functions, that is, functions $r$ that satisfy

$$0 < \liminf \frac{\log r(n)}{n} \leq \limsup \frac{\log r(n)}{n} < \infty,$$

it is known that (1.1) holds if and only if the Foster–Lyapunov drift condition holds that is, there exist an extended real-valued function $V : \mathsf{X} \to [1, \infty]$ finite at some $x_0 \in \mathsf{X}$, a petite set $C$, $\lambda \in (0,1)$, $b > 0$ and $c > 0$ such that $c^{-1} f \le V \le cf$ and

$$(1.2) \qquad\qquad PV \le \lambda V + b\mathbf{1}_C.$$

In that case, the convergence (1.1) holds for all $x$ in the set $\{V < \infty\}$ which is of $\pi$ measure 1. See, for instance, Meyn and Tweedie (1993) (hereafter, MT), Theorem 16.0.1.

For rates of convergence slower than geometric, no such definitive result exists. An important family of such rates is the class of subgeometric rate functions, defined in Nummelin and Tuominen (1983) as follows. Let $\Lambda_0$ be the set of positive nondecreasing functions $r_0$ such that $r_0(0) \ge 1$ and $\log\{r_0(n)\}/n$ decreases to 0. The class of subgeometric rate functions is the set $\Lambda$ of positive functions $r$ such that there exists a sequence $r_0 \in \Lambda_0$ that satisfies

$$(1.3) \qquad\qquad 0 < \liminf \frac{r(n)}{r_0(n)} \le \limsup \frac{r(n)}{r_0(n)} < \infty.$$

This class includes, for example, polynomial rate functions, that is, rate functions $r$ such that (1.3) holds with $r_0(n) := (1+n)^\beta$ for some $\beta > 0$. It also includes rate functions which increase faster than polynomially, for example, rate functions $r$ satisfying (1.3) with

$$(1.4) \qquad\qquad r_0(n) := \{1 + \log(n)\}^\alpha (n+1)^\beta e^{cn^\gamma}$$
$$\text{for } \alpha, \beta \in \mathbb{R}, \ \gamma \in (0,1) \text{ and } c > 0.$$

We will refer to these rates as subexponential in order to distinguish them in the broad class of subgeometric rates.

Tuominen and Tweedie (1994) [see also Nummelin and Tuominen (1983)] gave a set of equivalent conditions that imply the convergence (1.1) with a subgeometric rate function $r \in \Lambda$. To state this result, we first recall some notation and definitions.

A measurable set $C$ is $\psi_a$-petite (or petite) if there exist a distribution $a := (a(n), n \ge 0)$, a constant $\varepsilon > 0$ and a nontrivial measure $\psi_a$ on $\mathcal{B}(\mathsf{X})$ such that, for all $x \in C$, $B \in \mathcal{B}(\mathsf{X})$,

$$K_a(x,B) := \sum_{n \ge 0} a(n) P^n(x,B) \ge \psi_a(B).$$

The return time to a measurable set $A$ is defined as $\tau_A := \inf\{n \ge 1, \Phi_n \in A\}$ (with the convention $\inf \varnothing = +\infty$). Let $\psi$ be a maximal irreducibility measure and let $\mathcal{B}^+(\mathsf{X})$ be the class of accessible sets, that is, sets $B \in \mathcal{B}(\mathsf{X})$ such that $\psi(B) > 0$. A set $A \in \mathcal{B}(\mathsf{X})$ is called full if $\psi(A^c) = 0$, absorbing if $P(x,A) = 1$ for all $x \in A$ and, for a measurable positive function $f$ and a rate function $r$, $A$ is said to be $(f,r)$-regular if, for every $B \in \mathcal{B}^+(\mathsf{X})$,

$$\sup_{x \in A} \mathbb{E}_x \left[ \sum_{k=0}^{\tau_B - 1} r(k) f(\Phi_k) \right] < \infty.$$

A $(1, 1)$-regular set is simply said to be regular. A finite positive measure $\lambda$ on $\mathcal{B}(\mathsf{X})$ is said to be $(f, r)$-regular if $\mathbb{E}_\lambda[\sum_{k=0}^{\tau_B - 1} r(k) \, f(\Phi_k)] < \infty$ for all sets $B \in \mathcal{B}^+(\mathsf{X})$. If, for some $x \in \mathsf{X}$, the Dirak measure $\delta_x$ is $(f, r)$-regular, then the point $x$ is said to be $(f, r)$-regular. The set of all $(f, r)$-regular points is denoted by $S(f, r)$.

We can now recall (part of) Theorem 2.1 of Tuominen and Tweedie (1994).

THEOREM 1.1 [Tuominen and Tweedie (1994)]. *Assume that $P$ is $\psi$-irreducible and aperiodic. Let $f : \mathsf{X} \to [1, \infty]$ be a measurable function, and let $r \in \Lambda$ be given. The following conditions are equivalent.*

(i) *There exists a petite set $C \in \mathcal{B}(\mathsf{X})$ such that*

$$\sup_{x \in C} \mathbb{E}_x \left[ \sum_{k=0}^{\tau_C - 1} r(k) f(\Phi_k) \right] < \infty.$$

(ii) *There exist a sequence of extended real valued functions $(V_n, n \geq 0)$, $V_n : \mathsf{X} \to [1, \infty]$, a petite set $C \in \mathcal{B}(\mathsf{X})$ and a constant $b < \infty$ such that $V_0$ is bounded on $C$,*

$$V_0(x) = +\infty \quad \Longrightarrow \quad V_1(x) = +\infty$$

*and*

(1.5) $$PV_{n+1} + r(n)f \leq V_n + br(n)\mathbf{1}_C.$$

(iii) *There exists an $(f, r)$-regular set $A \in \mathcal{B}^+(\mathsf{X})$.*

*Any of these conditions implies that, for all $x \in S(f, r)$,*

$$r(n)\|P^n(x, \cdot) - \pi(\cdot)\|_f = 0, \qquad n \to \infty,$$

*and the set $S(f, r)$ is full, absorbing and contains the set $\{V_0 < \infty\}$. Moreover, for all $(f, r)$-regular initial distributions $\lambda$, $\mu$, there exists a constant $c$ such that*

$$\sum_{n=0}^\infty r(n) \iint \lambda(dx)\mu(dy)\|P^n(x, \cdot) - P^n(y, \cdot)\|_f$$
$$\leq c\big(\lambda(V_0) + \mu(V_0)\big).$$

This theorem cannot be improved since it provides a necessary and sufficient condition, but the sequence of drift conditions (1.5) is notoriously difficult to check in practice and one has very little insight on the way to choose the family of drift function $(V_n, n \geq 0)$. This is why these drift conditions, to the best of our knowledge, have seldom been used directly.

A first step toward finding a more practical drift condition was taken by Jarner and Roberts (2002), who, simplifying and generalizing an argument in Fort and Moulines (2000), have shown that if there exist a function $V : \mathsf{X} \to [1, \infty)$ finite at some $x_0 \in \mathsf{X}$, positive constants $b$ and $c$, a petite set $C$ and $\alpha \in [0, 1)$ such that

$$PV + cV^\alpha \leq V + b\mathbf{1}_C,$$

then the chain is positive recurrent and, for each $\beta \in [1, 1/(1-\alpha)]$, the convergence (1.1) holds for all $x \in \{V < \infty\}$ which is of $\pi$ measure 1, with $r(n) := n^{\beta-1}$ and $f := V^{1-\beta(1-\alpha)}$. It is noteworthy that there is a balance between the rate of convergence and the norm: the larger the norm, the slower the former. In particular, the fastest rate of convergence $[r(n) = n^{\alpha/(1-\alpha)}]$ corresponds to the total variation norm, and the slowest rate ($r \equiv 1$) corresponds to the $V^\alpha$ norm.

In this paper, we consider the following drift condition which generalizes the Foster–Lyapunov and the Jarner–Roberts drift conditions.

CONDITION $\mathbf{D}(\phi, V, C)$.   There exist a function $V : \mathsf{X} \to [1, \infty)$, a concave monotone nondecreasing differentiable function $\phi : [1, \infty) \mapsto (0, \infty]$, a measurable set $C$ and a finite constant $b$ such that

$$PV + \phi \circ V \leq V + b\mathbf{1}_C.$$

Here $\phi$ is assumed differentiable for convenience. It can be relaxed since a concave function has nonincreasing left and right derivatives everywhere. If $\mathbf{D}(\phi, V, C)$ holds for some petite set $C$ and there exists $x_0 \in \mathsf{X}$ such that $V(x_0) < \infty$, then the $f$-norm ergodic theorem (see MT, Theorem 14.0.1) states that there exists a unique invariant distribution $\pi$, $\pi(\phi \circ V) < \infty$ and

$$\lim_n \|P^n(x, \cdot) - \pi\|_{\phi \circ V} = 0,$$

for all $x$ in the set of $\pi$-measure 1 $\{V < \infty\}$. If, in addition, $\pi(V) < \infty$, then there exists a finite constant $B$ such that, for all $x \in \{V < \infty\}$,

$$\sum_{n=0}^{\infty} \|P^n(x, \cdot) - \pi\|_{\phi \circ V} \leq B(1 + V(x)).$$

The $(\phi \circ V)$-norm is the maximal norm for which convergence can be proved under Condition $\mathbf{D}(\phi, V, C)$ and in that case, the rate of convergence is minimal: $r \equiv 1$. This implies that, for any function $1 \leq f \leq \phi \circ V$, convergence in the $f$-norm also holds. In order to determine the rate of convergence in the $f$-norm by means of Theorem 1.1, we should try to find a sequence of functions $(V_n, n \geq 0)$ such that (1.5) holds, but this is precisely what we are trying to avoid doing for all functions $f$. Instead, having in mind the balance between the rate of convergence and the norm, we will first determine the rate of convergence in the total variation norm by using the criterion (1.5) and then interpolate intermediate rates of convergence in the $f$-norm.

Thus, this new drift condition not only generalizes former results, but also yields a very straightforward way of proving subgeometric rates of convergence, in particular subexponential rates of the form (1.4). The interpolation technique is a key tool to obtain these rates easily and in practice yields all usual rates.

The rest of the paper is organized as follows. Our main result, Theorem 2.8, is stated and proved in the next section. Several typical functions $\phi$ are then considered, leading to a variety of subgeometric rate functions. The flexibility of the interpolation technique is illustrated by exhibiting new pairs of rates and controls functions. Several applications are given in Section 3. We establish subgeometric rates of convergence (in particular, faster than polynomial rates) in several models: first-order nonlinear autoregressive models, stochastic unit root models and the random walk multidimensional Hastings–Metropolis algorithm, under conditions which do not imply geometric ergodicity. These examples should illustrate the efficiency of Condition $\mathbf{D}(\phi, V, C)$ and the fact that it is indeed easier to check than the sequence of drift conditions (1.5).

## 2. Main result.

2.1. *Rate of convergence in the total variation norm.* Let $\phi : [1, \infty) \to (0, \infty)$ be a concave nondecreasing differentiable function. Define

$$(2.1) \qquad H_\phi(v) := \int_1^v \frac{dx}{\phi(x)}.$$

Then $H_\phi$ is a nondecreasing concave differentiable function on $[1, \infty)$. Moreover, since $\phi$ is concave, $\phi'$ is nonincreasing. Hence $\phi(v) \le \phi(1) + \phi'(1)(v-1)$ for all $v \ge 1$, which implies that $H_\phi$ increases to infinity. We can thus define its inverse $H_\phi^{-1} : [0, \infty) \to [1, \infty)$, which is also an increasing and differentiable function, with derivative $(H_\phi^{-1})'(x) = \phi \circ H_\phi^{-1}(x)$. For $k \in \mathbb{N}$, $z \ge 0$ and $v \ge 1$, define

$$r_\phi(z) := (H_\phi^{-1})'(z) = \phi \circ H_\phi^{-1}(z),$$

$$(2.2) \qquad H_k(v) := \int_0^{H_\phi(v)} r_\phi(z+k) \, dz = H_\phi^{-1}(H_\phi(v) + k) - H_\phi^{-1}(k),$$

$$V_k := H_k \circ V.$$

We will show that, provided Condition $\mathbf{D}(\phi, C, V)$ holds with $C$ petite and $\sup_{x \in C} V(x) < \infty$, the chain $(\Phi_k, k \ge 0)$ is $(1, r_\phi)$-regular; that is, $r_\phi$ is the rate of convergence in total variation norm that can be deduced from the drift condition. To this end, we will use Theorem 1.1(ii); that is, we will show that (1.5) holds with $(V_k, k \ge 0)$, $f := 1$ and $r := r_\phi$.

PROPOSITION 2.1.  *Assume* $\mathbf{D}(\phi, V, C)$. *Then* $r_\phi$ *is log concave and, for all* $k \geq 0$, $H_k$ *is concave and*

$$PV_{k+1} \leq V_k - r_\phi(k) + \frac{br_\phi(k+1)}{r_\phi(0)}\mathbf{1}_C.$$

PROOF.  Note first that $r'_\phi(z)/r_\phi(z) = \phi' \circ H_\phi^{-1}(z)$ for all $z \geq 0$. Since $\phi'$ is nonincreasing and $H_\phi^{-1}$ is increasing, $\phi' \circ H_\phi^{-1}$ is nonincreasing and $\log(r_\phi)$ is concave. This implies that, for any fixed $k \geq 0$, the function $z \mapsto r_\phi(z + k)/r_\phi(z)$ is a decreasing function. The derivative of $H_k$ has the following expression:

$$\begin{aligned}
(2.3) \qquad H'_k(v) &= r_\phi(H_\phi(v) + k)/\phi(v) \\
&= r_\phi(H_\phi(v) + k)/r_\phi(H_\phi(v)).
\end{aligned}$$

Since $H_\phi$ is increasing, it follows from the discussion above that $H'_k$ is nonincreasing; hence $H_k$ is concave for all $k \geq 0$. Applying (2.3) and the fact that $r_\phi$ is increasing, we obtain

$$\begin{aligned}
&H_{k+1}(v) - H_k(v) \\
&= \int_0^{H_\phi(v)} \{r_\phi(z + k + 1) - r_\phi(z + k)\}\, dz \\
&= \int_0^{H_\phi(v)} \int_0^1 r'_\phi(z + k + s)\, ds\, dz \\
&= \int_0^1 \{r_\phi(H_\phi(v) + k + s) - r_\phi(k + s)\}\, ds \\
&\leq r_\phi(H_\phi(v) + k + 1) - r_\phi(k) \\
&= \phi(v)H'_{k+1}(v) - r_\phi(k).
\end{aligned}$$

We have thus shown the following inequality which is the key tool of the proof:

$$(2.4) \qquad H_{k+1}(v) - \phi(v)H'_{k+1}(v) \leq H_k(v) - r_\phi(k).$$

Let $g$ be a concave differentiable function on $[1, \infty)$. Since $g'$ is decreasing, for all $v \geq 1$ and $x \in \mathbb{R}$ such that $v + x \geq 1$, it holds that

$$(2.5) \qquad g(v + x) \leq g(v) + g'(v)x.$$

Applying this property to the concave function $H_{k+1}$, we obtain for all $k \geq 0$, $x \in \{V < \infty\}$,

$$\begin{aligned}
PV_{k+1}(x) &\leq H_{k+1}\{V(x) - \phi \circ V(x) + b\mathbf{1}_C(x)\} \\
&\leq H_{k+1}(V(x)) - \phi \circ V(x)H'_{k+1}(V(x)) + bH'_{k+1}(V(x))\mathbf{1}_C(x) \\
&\leq H_{k+1}(V(x)) - \phi \circ V(x)H'_{k+1}(V(x)) + bH'_{k+1}(1)\mathbf{1}_C(x).
\end{aligned}$$

Applying (2.3) and (2.4), we obtain that $H'_{k+1}(1) = r_\phi(k+1)/r_\phi(0)$ and

$$PV_{k+1}(x) \leq V_k(x) - r_\phi(k) + \frac{br_\phi(k+1)}{r_\phi(0)}\mathbf{1}_C(x).$$

This inequality still holds for $x \in \{V = \infty\}$, which concludes the proof. $\square$

The drift condition $\mathbf{D}(\phi, V, C)$ and Proposition 2.1 imply the following bounds for the modulated moments of the return time $\tau_C$, by application of Proposition 11.3.2 in MT.

PROPOSITION 2.2. *Assume* $\mathbf{D}(\phi, V, C)$. *Then, for all* $x \in \mathsf{X}$,

$$\mathbb{E}_x\left[\sum_{k=0}^{\tau_C-1} \phi \circ V(\Phi_k)\right] \leq V(x) + b\mathbf{1}_C(x),$$

$$\mathbb{E}_x\left[\sum_{k=0}^{\tau_C-1} r_\phi(k)\right] \leq V(x) + \frac{br_\phi(1)}{r_\phi(0)}\mathbf{1}_C(x).$$

In order to apply Theorem 1.1, we must also check the following conditions:

  (i) the rate sequence $r_\phi := (\phi \circ H_\phi^{-1}(k), k \geq 0)$ belongs to $\Lambda$,
  (ii) the drift function $V$ is bounded on $C$.

The next lemma gives a simple criterion to check that $r_\phi \in \Lambda$.

LEMMA 2.3. *If* $\lim_{t \to \infty} \phi'(t) = 0$, *then* $r_\phi \in \Lambda$.

PROOF. We have already noted that $r'_\phi(x)/r_\phi(x) = \phi' \circ H_\phi^{-1}(x)$ for all $x \geq 0$. Let $r$ be any differentiable function such that $r(0) = 1$ and $\lim_{x \to \infty} r'(x)/r(x) = 0$. Then, applying Cesaro's lemma, we obtain,

$$\frac{\log(r(n))}{n} = \frac{1}{n}\int_0^n \frac{r'(s)}{r(s)}\,ds \to 0.$$

If, moreover, $r'/r$ decreases, then $\log(r(x))/x$ also decreases. Thus $r_\phi \in \Lambda$. $\square$

If the condition $\sup_{x \in C} V(x) < \infty$ is not satisfied and if the set $C$ is petite, the drift condition $\mathbf{D}(\phi, V, C)$ can be slightly modified so that it holds with a new set $C$ on which $V$ is bounded. The following lemma, adapted from Theorem 14.2.3 of MT, states this formally.

LEMMA 2.4. *Assume that* $\mathbf{D}(\phi, V, C)$ *holds for some petite set* $C$ *and that* $\lim_{v \to \infty} \phi(v) = \infty$. *Then, for all* $M \geq 1$, *the sublevel sets* $\{x \in \mathsf{X}, V(x) \leq M\}$ *are petite. In addition, for any* $\beta$, $0 < \beta < 1$, *there exists a sublevel set* $C_\beta$ *such that* $\mathbf{D}(\beta\phi, V, C_\beta)$ *holds.*

PROOF.    Since $\phi$ is positive nondecreasing and $V \geq 1$, Condition $\mathbf{D}(\phi, V, C)$ implies the drift condition $PV \leq V - \phi(1) + b\mathbf{1}_C$. Theorem 11.3.11 of MT shows that, for all accessible sets $B \in \mathcal{B}^+(\mathsf{X})$, there exists a constant $c(B) < \infty$ such that, for all $x \in \mathsf{X}$, we have $\phi(1)\mathbb{E}_x[\tau_B] \leq V(x) + c(B)$. Hence, every set $A \in \mathcal{B}(\mathsf{X})$ such that $\sup_{x \in A} V(x) < \infty$ is $(1, 1)$-regular, and the sublevel sets are all $(1, 1)$-regular. Proposition 11.3.8 of MT shows that if a set $A$ is $(1, 1)$-regular, then it is petite. Hence, all the sublevel sets are petite.

Since $\lim_{v \to \infty} \phi(v) = \infty$ for all $\beta \in (0, 1)$, there exists $M_\beta$ such that $v > M_\beta$ implies $\phi(v) \geq b/(1 - \beta)$. For $x \notin C_\beta := \{V \leq M_\beta\}$, we thus have $b \leq (1 - \beta)\phi(V(x))$ and

$$PV + \beta\phi(V) \leq V + (\beta - 1)\phi(V) + b\mathbf{1}_C \leq V.$$

For $x \in C_\beta$, since $\beta \in (0, 1)$, it trivially holds that

$$PV + \beta\phi(V) \leq V + b. \qquad \square$$

We are now in position to establish the rate of convergence in total variation distance.

PROPOSITION 2.5.    *Let $P$ be a $\psi$-irreducible and aperiodic kernel. Assume that $\mathbf{D}(\phi, V, C)$ holds for a function $\phi$ such that $\lim_{t \to \infty} \phi'(t) = 0$, a petite set $C$ and a function $V$ such that $\{V < \infty\} \neq \varnothing$. Then, there exists an invariant probability measure $\pi$, and for all $x$ in the full and absorbing set $\{V < \infty\}$,*

$$\lim_n r_\phi(n)\|P^n(x, \cdot) - \pi(\cdot)\|_{TV} = 0.$$

*Any probability measure $\lambda$ such that $\lambda(V) < \infty$ is $(1, r_\phi)$-regular and for two $(1, r_\phi)$-regular distributions $\lambda, \mu$, there exists a constant $c$ such that*

$$\sum_{n=0}^{\infty} r_\phi(n) \iint \lambda(dx)\mu(dy)\|P^n(x, \cdot) - P^n(y, \cdot)\|_{TV}$$
$$\leq c\big(\lambda(V) + \mu(V)\big).$$

REMARK 1.    Since $\phi'$ is nonincreasing, if we do not assume that $\lim_{v \to \infty} \phi'(v) = 0$, then there exists $c \in (0, 1)$ such that $\lim_{v \to \infty} \phi'(v) = c > 0$. This yields $v - \phi(v) \leq (1 - c)v + c - \phi(1)$. In this case, Condition $\mathbf{D}(\phi, V, C)$ implies the Foster–Lyapunov drift condition, and the chain is geometrically ergodic.

PROOF OF PROPOSITION 2.5.    The only statement which requires a proof is the fact that any probability measure such that $\lambda(V) < \infty$ is $(1, r_\phi)$-regular. This assertion is established in Proposition 3.1(ii) of Tuominen and Tweedie (1994), and relies on Lemma 3.1 of Nummelin and Tuominen (1983). We nevertheless propose a proof that shortens the previous one. The proof is adapted from the proof of Theorem 14.2.3 of MT. Proposition 2.1 shows that there exist a sequence of drift functions $(V_k, k \geq 0)$ and a constant $b$ such that $V_0 \leq V$ and

$$PV_{k+1} \leq V_k - r_\phi(k) + b\phi(1)^{-1} r_\phi(k+1) \mathbf{1}_C.$$

Dynkin's formula shows that, for all accessible set $B$,

$$\mathbb{E}_x \left[ \sum_{k=0}^{\tau_B - 1} r_\phi(k) \right]$$

$$\leq V_0(x) + b\phi(1)^{-1} \mathbb{E}_x \left[ \sum_{k=0}^{\tau_B - 1} r_\phi(k+1) \mathbf{1}_C(\mathbf{\Phi}_k) \right].$$

From Propositions 5.5.5 and 5.5.6 of MT, we can assume without loss of generality that $C$ is $\psi_a$-petite, where $\psi_a$ is equivalent to $\psi$, and that the sampling distribution $a$ has finite mean $m_a := \sum_{j=1}^{\infty} j a_j < \infty$. By the Comparison Theorem (MT, Theorem 14.2.2), the bound $\mathbf{1}_C(x) \leq \psi_a(B)^{-1} K_a(x, B)$ and the fact that $r_\phi$ is nondecreasing, we have

$$\mathbb{E}_x \left[ \sum_{k=0}^{\tau_B - 1} r_\phi(k) \right]$$

$$\leq V_0(x) + b\phi(1)^{-1} \mathbb{E}_x \left[ \sum_{k=0}^{\tau_B - 1} r_\phi(k+1) \mathbf{1}_C(\mathbf{\Phi}_k) \right]$$

$$\leq V_0(x) + b\phi(1)^{-1} \psi_a(B)^{-1} \sum_{i \geq 0} a_i \mathbb{E}_x \left[ \sum_{k=0}^{\tau_B - 1} r_\phi(k+1) \mathbf{1}_B(\mathbf{\Phi}_{k+i}) \right]$$

$$\leq V_0(x) + b\phi(1)^{-1} \psi_a(B)^{-1} m_a \mathbb{E}_x[r_\phi(\tau_B)].$$

For $k \geq 1$, define $R_\phi(k) := \sum_{j=0}^{k-1} r_\phi(j)$. Since $r_\phi$ is subgeometric, it holds that $\lim_{k \to \infty} r_\phi(k)/R_\phi(k) = 0$. Hence, for any $\delta > 0$, there exists a constant $c(\delta)$ such that, for all $k \geq 1$, $r_\phi(k) \leq \delta R_\phi(k) + c(\delta)$. This yields

$$\mathbb{E}_x[R_\phi(\tau_B)] \leq V_0(x) + b\phi(1)^{-1} \psi_a(B)^{-1} m_a \big( \delta \mathbb{E}_x[R_\phi(\tau_B)] + c(\delta) \big).$$

Thus for small enough $\delta$, we obtain

$$(2.6) \qquad \mathbb{E}_x[R_\phi(\tau_B)] \leq \frac{V_0(x) + bm_a \psi_a^{-1}(B)c(\delta)\phi(1)^{-1}}{1 - b\delta m_a \psi_a^{-1}(B)\phi(1)^{-1}}. \qquad \square$$

2.2. *Rate of convergence in $f$-norms.* As already mentioned in the polynomial case and discussed in Tuominen and Tweedie (1994), in the subgeometric case there is a compromise between the rate of convergence and the control function. In what follows, we will show that it is possible at almost no cost to obtain many intermediate rates of convergence and control functions. Let $\mathcal{Y}$ be the set of pairs of ultimately nondecreasing functions $\Psi_1$ and $\Psi_2$ defined on $[1, \infty)$ such that $\lim_{x \to \infty} \Psi_1(x) = \infty$ or $\lim_{x \to \infty} \Psi_2(x) = \infty$ and, for all $x, y \in [1, \infty)$,

$$(2.7) \qquad \Psi_1(x)\Psi_2(y) \le x + y.$$

The set $\mathcal{Y}$ includes, for example, $\Psi_1(x) = x$ and $\Psi_2(x) = 1$, but there are of course more interesting examples. For example, it is well known that, for any $x, y \ge 0$, and $p$ and $q$ such that $1/p + 1/q = 1$, we have

$$xy \le x^p/p + y^q/q.$$

Hence, the pair of functions $\Psi_1(x) = p^{1/p}x^{1/p}$, $\Psi_2(x) = q^{1/q}x^{1/q}$ satisfies (2.7). These are precisely the interpolating functions used in Jarner and Roberts (2002) to derive polynomial rates of convergence. Young functions provide many useful interpolating functions. We recall their definition. Let $\varrho_1 : (0, \infty) \to (0, \infty)$ be an increasing left-continuous function such that $\varrho_1(0) = 0$ and $\lim_{v \to +\infty} \varrho_1(v) = +\infty$. Let $\varrho_2$ be the left-continuous inverse of $\varrho_1$, which is increasing and satisfies also $\varrho_2(0) = 0$ and $\lim_{v \to +\infty} \varrho_2(v) = +\infty$. Define then $G_i(x) := \int_0^x \varrho_i(t)\, dt$, $i = 1, 2$. The well-known Young inequality states that, for all $x, y \ge 0$, we have

$$(2.8) \qquad xy \le G_1(x) + G_2(y).$$

Let $\Psi_i$ be the inverse of $G_i$, $i = 1, 2$. Then $\Psi_1$ and $\Psi_2$ are concave functions and it follows immediately from (2.8) that the pair $(\Psi_1, \Psi_2)$ satisfies (2.7).

We use this full scale of interpolating functions in combination with Proposition 2.2 to derive bounds for the modulated moment of return time to the set $C$. More precisely, we have the following.

PROPOSITION 2.6.  *Assume* $\mathbf{D}(\phi, V, C)$ *and let* $(\Psi_1, \Psi_2) \in \mathcal{Y}$. *Then*

$$\mathbb{E}_x\left[\sum_{k=0}^{\tau_C - 1} \Psi_1\big(r_\phi(k)\big)\Psi_2\big(\phi \circ V(\Phi_k)\big)\right]$$
$$\le 2V(x) + b\big(1 + r_\phi(1)/r_\phi(0)\big)\mathbf{1}_C(x).$$

We need a criterion for a rate function $\Psi_1 \circ r_\phi$ to be subgeometric. Note that if the pair $(\Psi_1, \Psi_2)$ belongs to $\mathcal{Y}$, then, for large enough $x$, it holds that $\Psi_i(x) \le 2x$, $i = 1, 2$.

LEMMA 2.7. *Assume that* $\lim_{t\to\infty}\phi'(t)=0$. *For any nondecreasing function* $\Psi$ *such that* $\Psi(x)\le ax$ *for some constant* $a$, $\Psi\circ r_\phi\in\Lambda$.

The next theorem summarizes all our previous results.

THEOREM 2.8. *Let* $P$ *be a* $\psi$-*irreducible and aperiodic kernel. Assume that* $\mathbf{D}(\phi,V,C)$ *holds for a function* $\phi$ *such that* $\lim_{t\to\infty}\phi'(t)=0$ *and a petite set* $C$ *such that* $\sup_C V<\infty$. *Let* $(\Psi_1,\Psi_2)\in\mathcal{Y}$. *Then, there exists an invariant probability measure* $\pi$, *and for all* $x$ *in the full set* $\{V<\infty\}$,

$$\lim_n\Psi_1(r_\phi(n))\|P^n(x,\cdot)-\pi(\cdot)\|_{\Psi_2(\phi\circ V)}=0.$$

*Any probability measure* $\lambda$ *such that* $\lambda(V)<\infty$ *is* $(\Psi_2(\phi\circ V),\Psi_1(r_\phi))$-*regular, and for two such distributions* $\lambda,\mu$, *there exists a constant* $c$ *such that*

$$\sum_{n=0}^\infty\Psi_1(r_\phi(n))\iint\lambda(dx)\mu(dy)\|P^n(x,\cdot)-P^n(y,\cdot)\|_{\Psi_2(\phi\circ V)}$$
$$\le c(\lambda(V)+\mu(V)).$$

PROOF. From Proposition 2.6 we have

$$\sup_{x\in C}\mathbb{E}_x\left[\sum_{k=0}^{\tau_C-1}\Psi_1(r_\phi(k))\Psi_2(\phi\circ V(\Phi_k))\right]<\infty.$$

Theorem 1.1 shows that $\Phi$ is $(\Psi_2(\phi\circ V),\Psi_1(r_\phi))$-regular. As in the proof of Proposition 2.5, and using again the comparison theorem, for any set $B\in\mathcal{B}^+(\mathsf{X})$, there exist constants $c_1(B)$ and $c_2(B)$ such that

$$\mathbb{E}_x\left[\sum_{k=0}^{\tau_B-1}\phi\circ V(\Phi_k)\right]+\mathbb{E}_x\left[\sum_{k=0}^{\tau_B-1}r_\phi(k)\right]\le c_1(B)V(x)+c_2(B).$$

Hence, for any $(\Psi_1,\Psi_2)\in\mathcal{Y}$, we have

$$\mathbb{E}_x\left[\sum_{k=0}^{\tau_B-1}\Psi_1(r_\phi(k))\Psi_2(\phi\circ V(\Phi_k))\right]\le c_1(B)V(x)+c_2(B),$$

which shows that any probability measure such that $\lambda(V)<\infty$ is $(\Psi_2(\phi\circ V),\Psi_1(r_\phi))$-regular. □

2.3. *Some usual rate functions.* In this section, we provide examples of rates of convergence obtained by Theorem 2.8 for several functions $\phi$. In Section 3, we will provide examples and explicitly determine the drift function $V$ and the set $C$. For two sequences $u_n$ and $v_n$, we write $u_n\asymp v_n$, if there exist positive constants $c_1$ and $c_2$ such that, for large $n$, $c_1u_n\le v_n\le c_2u_n$.

We assume throughout this section that Condition $\mathbf{D}(\phi,V,C)$ holds for some petite set $C$ and $\sup_C V<\infty$.

*Polynomial rates of convergence.* Polynomial rates of convergence have been widely studied recently under various conditions [see Veretennikov (1997, 2000), Tanikawa (2001), Jarner and Roberts (2002) and Fort and Moulines (2003)]. As already mentioned, polynomial rates of convergence are associated to the functions $\phi(v) := cv^\alpha$ for some $\alpha \in [0, 1)$ and $c \in (0, 1]$ and the rate of convergence in total variation distance is $r_\phi(n) \propto n^{\alpha/(1-\alpha)}$. Set $\Psi_1(x) := ((1-p)x)^{(1-p)}$ and $\Psi_2(x) := (px)^p$ for some $p$, $0 < p < 1$. Applying Theorem 2.8 yields, for $x \in \{V < \infty\}$,

$$(2.9) \qquad \lim_n n^{(1-p)\alpha/(1-\alpha)} \|P^n(x, \cdot) - \pi(\cdot)\|_{V^{\alpha p}} = 0.$$

This convergence remains valid for $p = 0, 1$ by Proposition 2.2. Set $\kappa := 1 + (1-p)\alpha/(1-\alpha)$ so that $1 \le \kappa \le 1/(1-\alpha)$. With this notation, (2.9) reads

$$\lim_n n^{\kappa-1} \|P^n(x, \cdot) - \pi(\cdot)\|_{V^{1-\kappa(1-\alpha)}} = 0,$$

which is Theorem 3.6 of Jarner and Roberts (2002).

It is possible to extend this result by using more general interpolation functions. We can, for example, obtain nonpolynomial rates of convergence with control functions which are not simply power of the drift functions. To illustrate this point, set for $b > 0$, $\Psi_1(x) := (1 \vee \log(x))^b$ and $\Psi_2(x) := x(1 \vee \log(x))^{-b}$. It is not difficult to check that we have

$$\sup_{(x,y) \in [1,\infty) \times [1,\infty)} (x+y)^{-1} \Psi_1(x) \Psi_2(y) < \infty,$$

so that, for all $x \in \{V < \infty\}$, we have

$$(2.10) \qquad \lim_n \log^b(n) \|P^n(x, \cdot) - \pi(\cdot)\|_{V^\alpha(1+\log(V))^{-b}} = 0,$$

$$(2.11) \qquad \lim_n n^{\alpha/(1-\alpha)} \log^{-b}(n) \|P^n(x, \cdot) - \pi(\cdot)\|_{(1+\log(V))^b} = 0,$$

and for all $0 < p < 1$,

$$\lim_n n^{(1-p)\alpha/(1-\alpha)} \log^b n \|P^n(x, \cdot) - \pi(\cdot)\|_{V^{\alpha p}(1+\log V)^{-b}} = 0.$$

*Logarithmic rates of convergence.* We now consider drift conditions which imply rates of convergence slower than any polynomial. Such rates are obtained when Condition $\mathbf{D}(\phi, V, C)$ holds with a function $\phi$ that increases to infinity slower than polynomially. We only consider here the case $\phi(v) = c(1 + \log(v))^\alpha$ for some $\alpha \ge 0$ and $c \in (0, 1]$. A straightforward calculation shows that

$$r_\phi(n) \asymp \log^\alpha(n).$$

Proposition 2.5 shows that the chain is $(1, \log^\alpha(n))$ and $((1 + \log V)^\alpha, 1)$-regular. Applying Theorem 2.8, intermediate rates can be obtained along the same lines as above. Choosing, for instance, $\Psi_1(x) := ((1-p)x)^{1-p}$ and $\Psi_2(x) := (px)^p$ for $0 \leq p \leq 1$, the chain is $((1 + \log V)^{p\alpha}, \log(n)^{(1-p)\alpha})$-regular and thus, for all $x \in \{V < \infty\}$,

$$\lim_{n \to \infty} (1 + \log(n))^{(1-p)\alpha} \| P^n(x, \cdot) - \pi(\cdot) \|_{(1+\log(V))^{p\alpha}} = 0.$$

*Subexponential rates of convergence.* Subexponential rates [as defined in (1.4)] have been considered only recently in the literature. An example (in continuous time) has been studied by Malyshkin (2001); discrete-time examples are considered in the recent work by Klokov and Veretennikov (2002). These rates, which increase to infinity faster than polynomially, are obtained when Condition $\mathbf{D}(\phi, V, C)$ holds with $\phi$ such that, $v/\phi(v)$ goes to infinity slower than polynomially. More precisely, assume that $\phi$ is concave and differentiable on $[1, +\infty)$ and that for large $v$, $\phi(v) = cv/\log^\alpha(v)$ for some $\alpha > 0$ and $c > 0$. A simple calculation yields

$$r_\phi(n) \asymp n^{-\alpha/(1+\alpha)} \exp(\{c(1+\alpha)n\}^{1/(1+\alpha)}),$$

and thus the chain is $(1, n^{-\alpha/(1+\alpha)} \exp(\{c(1+\alpha)n\}^{1/(1+\alpha)}))$ and $(V/(1+\log V)^\alpha, 1)$-regular. Applying Theorem 2.8 with $\Psi_1(x) := x^{1-p}(1 \vee \log(x))^{-b}$ and $\Psi_2(x) := x^p(1 \vee \log(x))^b$ for $p \in (0, 1)$ and $b \in \mathbb{R}$; $p = 0$ and $b > 0$; or $p = 1$ and $b < -\alpha$ yields, for all $x \in \{V < \infty\}$,

$$\begin{aligned}
(2.12) \quad & \lim_n n^{-(\alpha+b)/(1+\alpha)} \exp((1-p)\{c(1+\alpha)n\}^{1/(1+\alpha)}) \\
& \quad \times \| P^n(x, \cdot) - \pi(\cdot) \|_{V^p(1+\log V)^b} = 0.
\end{aligned}$$

Asymptotically, the term $n^{-(\alpha+b)/(1+\alpha)}$ is not very important and we can express (2.12) in a simpler way: for all $d < (1-p)\{c(1+\alpha)\}^{1/(1+\alpha)}$, we have, for the same values of $p$ and $b$,

$$\lim_n e^{dn^{1/(1+\alpha)}} \| P^n(x, \cdot) - \pi(\cdot) \|_{V^p(1+\log V)^b} = 0.$$

## 3. Applications.

We now illustrate our findings by applying Theorem 2.8 to several models. In Section 3.1, we exhibit a simple example where the rates obtained in Theorem 2.8 can be proved optimal. In the next sections we study several examples where no such optimality results are available.

In this section, $\langle \cdot, \cdot \rangle$ and $| \cdot |$ denote, respectively, the scalar product and the Euclidean norm in any Euclidean space. The transpose of a vector $v$ is denoted $v'$. If $u$ is a twice continuously differentiable real valued function on $\mathbb{R}^d$, $\nabla u$ (resp. $\nabla^2 u$) denotes its gradient (resp. its Hessian matrix).

3.1. *Backward recurrence time chain.* The backward recurrence time chain (see MT, Section 3.3.1) is a rich source of simple examples of stable and unstable behavior. We consider it here to provide examples of chains satisfying Condition $\mathbf{D}(\phi, V, C)$ and for which the rates of convergence implied by it are optimal.

Let $(p_n, n \geq 0)$ be a sequence of positive real numbers such that $p_0 = 1$, $p_n \in (0, 1)$ for all $n \geq 1$ and $\lim_{n \to \infty} \prod_{i=1}^n p_i = 0$. Consider the backward recurrence time chain $\mathbf{\Phi}$ with transition kernel $P$ defined as $P(n, n+1) = 1 - P(n, 0) = p_n$, for all $n \geq 0$. Then $\mathbf{\Phi}$ is irreducible and strongly aperiodic and $\{0\}$ is an atom. Let $\tau_0$ be the return time to $\{0\}$. We have, for all $n \geq 1$,

$$\mathbb{P}_0(\tau_0 = n + 1) = (1 - p_n) \prod_{j=0}^{n-1} p_j \quad \text{and} \quad \mathbb{P}_0(\tau_0 > n) = \prod_{j=0}^{n-1} p_j.$$

By Kac's theorem (MT, Theorem 10.2.2), since $\mathbf{\Phi}$ is $\psi$-irreducible and aperiodic, $\mathbf{\Phi}$ is positive recurrent if and only if $\mathbb{E}_0[\tau_0] < \infty$, that is,

$$\sum_{n=1}^\infty \prod_{j=1}^n p_j < \infty,$$

and the stationary distribution $\pi$ is given, by $\pi(0) = \pi(1) = 1/\mathbb{E}_0[\tau_0]$ and for $j \geq 2$,

$$\pi(j) = \frac{\mathbb{E}_0[\sum_{k=1}^{\tau_0} \mathbf{1}_{\{\Phi_k = j\}}]}{\mathbb{E}_0[\tau_0]} = \frac{\mathbb{P}_0(\tau_0 \geq j)}{\mathbb{E}_0[\tau_0]} = \frac{p_0 \cdots p_{j-2}}{\sum_{n=1}^\infty p_1 \cdots p_n}.$$

Because the distribution of the return time to the atom $\{0\}$ has such a simple expression in terms of the transition probability $(p_n, n \geq 0)$, we are able to exhibit the largest possible rate function $r$ such that the $(1, r)$-modulated moment of the return time $\mathbb{E}_0[\sum_{k=0}^{\tau_0-1} r(k)]$ is finite. We will also prove that the drift condition $\mathbf{D}(\phi, V, C)$ holds for appropriately chosen functions $V$ and $\phi$ and yields the optimal rate of convergence. Note also that, for any function $f$, it holds that

$$\mathbb{E}_0\left[\sum_{k=0}^{\tau_0-1} f(\Phi_k)\right] = \mathbb{E}_0\left[\sum_{k=0}^{\tau_0-1} f(k)\right].$$

Therefore there is no loss of generality to consider only $(1, r)$-modulated moments of the return time to zero.

If $\sup_{n \geq 1} p_n \leq \lambda < 1$, then, for all $\lambda < \mu < 1$, $\mathbb{E}_0[\mu^{-\tau_0}] < \infty$ and $\{0\}$ is thus a geometrically ergodic atom (MT, Theorem 15.1.5). Subgeometric rates of convergence in total variation norm are obtained when $\limsup p_n = 1$. Depending on the rate at which $p_n$ approaches 1, different behaviors can be obtained, covering essentially the three typical rates (polynomial, logarithmic and subexponential) discussed above.

*Polynomial rates.* Assume first that, for $\theta > 0$ and large $n$, $p_n = 1 - (1+\theta)n^{-1}$. Then $\prod_{i=1}^{n} p_i \asymp n^{-1-\theta}$. Thus, $\mathbb{E}_0[\sum_{k=0}^{\tau_0 - 1} r(k)] < \infty$ if and only if $\sum_{k=1}^{\infty} r(k)k^{-1-\theta} < \infty$. For instance, $r(n) := n^{\beta}$ with $0 \leq \beta < \theta$ is suitable.

*Logarithmic rates.* If for $\theta > 0$ and large $n$, $p_n = 1 - 1/n - (1+\theta)/(n\log(n))$, then $\prod_{j=1}^{n} p_j \asymp n^{-1}\log^{-1-\theta}(n)$, which is a summable series. Hence if $r$ is non-decreasing and $\sum_{k=1}^{\infty} r(k) \prod_{j=1}^{n} p_j < \infty$, then $r(k) = o(\log^{\theta}(k))$. In particular, $r(k) := \log^{\beta}(k)$ is suitable for all $0 \leq \beta < \theta$.

*Subgeometric rates.* If for large $n$, $p_n = 1 - \theta\beta n^{\beta-1}$ for some $\theta > 0$ and $\beta \in (0,1)$, then $\prod_{i=1}^{n} p_i \asymp e^{-\theta n^{\beta}}$. Thus, $\mathbb{E}_0[\sum_{k=0}^{\tau_0 - 1} e^{ak^{\beta}}] < \infty$ if $a < \theta$, and $\mathbb{E}_0[\sum_{k=0}^{\tau_0 - 1} e^{ak^{\beta}}] = \infty$ if $a \geq \theta$.

*Checking Condition $\mathbf{D}(\phi, V, C)$.* In order to prove that Proposition 2.5 provides the optimal rates of convergence, we now compute in each of the previous examples the rates of convergence it yields.

For the polynomial and subexponential cases, the same technique can be used. For $\gamma \in (0,1)$ and $x \in \mathbb{N}^*$, define $V(0) := 1$ and $V(x) := \prod_{j=0}^{x-1} p_j^{-\gamma}$. Then, for all $x \geq 0$, we have,

$$PV(x) = p_x V(x+1) + (1-p_x)V(0)$$
$$= p_x^{1-\gamma} V(x) + 1 - p_x$$
$$\leq V(x) - (1 - p_x^{1-\gamma})V(x) + 1 - p_x.$$

Thus, for $0 < \delta < 1 - \gamma$ and large enough $x$, it holds that

$$PV(x) \leq V(x) - \delta(1 - p_x)V(x). \tag{3.1}$$

1. Case $p_n = 1 - (1+\theta)n^{-1}$, $\theta > 0$. Then $V(x) \asymp x^{\gamma(1+\theta)}$ and $(1-p_x)V(x) \asymp V(x)^{1-1/(\gamma(1+\theta))}$. Thus Condition $\mathbf{D}(\phi, V, C)$ holds with $\phi(v) = cv^{\alpha}$ for $\alpha = 1 - 1/(\gamma(1+\theta))$ for any $\gamma \in (0,1)$. Theorem 2.8 yields the rate of convergence $n^{\alpha/(1-\alpha)} = n^{\gamma(1+\theta)-1}$, that is, $n^{\beta}$ for any $0 \leq \beta < \theta$.

2. Case $p_n = 1 - \theta\beta n^{\beta-1}$. Then, for large enough $x$, (3.1) yields

$$PV(x) \leq V(x) - \theta\beta\delta x^{\beta-1}V(x)$$
$$\leq cV(x)\{\log(V(x))\}^{1-1/\beta},$$

for $c < \theta^{1/\beta}\beta\delta$. Defining $\alpha := 1/\beta - 1$, Proposition 2.1 yields the following rate of convergence in total variation norm:

$$n^{-\alpha/(1+\alpha)}\exp(\{c(1+\alpha)n\}^{1/(1+\alpha)}) = n^{\beta-1}\exp(\theta\delta^{\beta}n^{\beta}).$$

Since $\delta$ is arbitrarily close to 1, we recover the fact that $\mathbb{E}_0[\sum_{k=0}^{\tau_0 - 1} e^{ak^{\beta}}] < \infty$ for any $a < \theta$.

3. Case $p_n = 1 - n^{-1} - (1+\theta)n^{-1}\log^{-1}(n)$, $\theta > 0$. Choose $V(x) := (\prod_{j=0}^{x-1} p_j)/\log^\varepsilon(x)$ for $\varepsilon > 0$ arbitrarily small. Then, for constants $c < c' < c'' < 1$ and large $x$, we obtain

$$PV(x) = \frac{\log^\varepsilon(x)}{\log^\varepsilon(x+1)}V(x) + 1 - p_x$$

$$= V(x) - c''\varepsilon\frac{V(x)}{x\log(x)} + 1 - p_x$$

$$\leq V(x) - c'\varepsilon\log^{\theta-\varepsilon}(x)$$

$$\leq V(x) - c\varepsilon\log^{\theta-\varepsilon}(V(x)).$$

Here again Theorem 2.8 yields the optimal rate of convergence.

3.2. *Symmetric random walk Hastings–Metropolis algorithm.* We consider the symmetric random walk Hastings–Metropolis algorithm. The purpose of this algorithm is to simulate from a probability distribution $\pi$ which is known only up to a scale factor. At each iteration, a move is proposed according to a random walk whose increment distribution has a symmetric density $q$ with respect to the Lebesgue measure $\mu_d$ on $\mathbb{R}^d$. The move is accepted with probability $\alpha(x, y)$ defined by

$$(3.2) \qquad \alpha(x, y) := \begin{cases} \min\left\{\dfrac{\pi(y)}{\pi(x)}, 1\right\}, & \text{if } \pi(x) > 0, \\ 1, & \text{if } \pi(x) = 0. \end{cases}$$

The transition kernel of the Metropolis algorithm is then given by

$$P(x, A) = \int_A \alpha(x, x+y)q(y)\,d\mu_d(y)$$

$$+ \mathbf{1}_A(x)\int\bigl(1 - \alpha(x, x+y)\bigr)q(y)\,d\mu_d(y).$$

Mengersen and Tweedie (1996) have shown that a real valued Metropolis chain is geometrically ergodic when the proposal density $q$ satisfies moment conditions and the target density $\pi$ is continuous, positive and log concave in the tails. This condition is necessary in the sense that if the chain is geometrically ergodic, then $\int \exp(s|z|)\pi(z)\,d\mu_d(z) < \infty$ for some $s > 0$. These results have been extended to the multidimensional case by Roberts and Tweedie (1996) and Jarner and Hansen (2000). Polynomial ergodicity was proved by Fort and Moulines (2000) for a target density with regularly varying tails. We now state conditions that imply subexponential rates of convergence.

ASSUMPTION 3.1. The target density $\pi$ is continuous and positive on $\mathbb{R}^d$ and there exist $m \in (0, 1)$, $r \in (0, 1)$, positive constants $d_i, D_i, i = 0, 1, 2$ and

$R_0 < \infty$ such that, if $|x| \geq R_0$, $x \mapsto \pi(x)$ is twice continuously differentiable and

$$(3.3) \qquad \left\langle \frac{\nabla \pi(x)}{|\nabla \pi(x)|}, \frac{x}{|x|} \right\rangle \leq -r,$$

$$(3.4) \qquad d_0|x|^m \leq -\log \pi(x) \leq D_0|x|^m,$$

$$(3.5) \qquad d_1|x|^{m-1} \leq |\nabla \log \pi(x)| \leq D_1|x|^{m-1},$$

$$(3.6) \qquad d_2|x|^{m-2} \leq |\nabla^2 \log \pi(x)| \leq D_2|x|^{m-2}.$$

The Weibull distribution on $\mathbb{R}$ with density $\pi(x) := \beta\gamma x^{\gamma-1} \exp(-\beta x^\gamma)$, for $x > 0$, $\beta > 0$ and $0 < \gamma < 1$ satisfies Assumption 3.1. Multidimensional examples are provided in Fort and Moulines (2000). For the sake of simplicity, we make the following assumption on the proposal density $q$.

ASSUMPTION 3.2. The proposal density $q$ is symmetric and bounded away from zero in a neighborhood of zero and is compactly supported; that is, there exists $c(q)$ such that, for all $|y| \geq c(q)$, $q(y) = 0$.

THEOREM 3.1. *Under Assumptions 3.1 and 3.2, there exist $z > 0$, $c > 0$, $r > 0$ such that $\mathbf{D}(\phi, V, C)$ holds with $V := \pi^{-z}$ and $\phi(v) := cv(1 + \log v)^{-2(1-m)/m}$ and $C := \{|x| \leq r\}$.*

Under Assumptions 3.1 and 3.2, Theorem 2.2 of Roberts and Tweedie (1996) shows that the chain is $\psi$-irreducible and aperiodic and nonempty bounded sets of $\mathcal{B}^+(\mathbb{R}^d)$ are petite. Thus, we obtain the following corollary.

COROLLARY 3.2. *There exist $c > 0$ and $z > 0$ such that, any probability measure $\lambda$ on $\mathbb{R}^d$ satisfying $\lambda(V) < \infty$ is $(f, r)$-regular with $r(n) = e^{cn^{m/(2-m)}}$ and $f = \pi^{-z}$.*

REMARK 2. Our result complements the work of Fort and Moulines (2000) who show that under Assumptions 3.1 and 3.2, the chain $\Phi$ is $(f, r)$-ergodic with $f(x) := (1 + |x|^\mu)$ and $r(n) := (1 + n)^\nu$, for any $\mu > 0$ and $\nu \geq 0$.

REMARK 3. The compactness Assumption 3.2 can probably be relaxed and replaced by an appropriate moment condition.

REMARK 4. We do not provide explicit values of the constants $c$ and $z$ here; these values can be deduced explicitly from the proof. It should be stressed that optimal values of these constants are related: the larger $c$, the smaller $z$ and vice versa. The same comments apply to Corollaries 3.4 and 3.6.

PROOF OF THEOREM 3.1.    Define $\mathcal{R}(x) := \{y \in \mathbb{R}^d, \pi(x+y) \le \pi(x)\}$ the potential rejection region. Using the definition of the transition kernel $P$, we have

$$PV(x) - V(x)$$
$$= \int \big(V(x+y) - V(x)\big) q(y) \, d\mu_d(y)$$
$$+ \int_{\mathcal{R}(x)} \big(V(x+y) - V(x)\big) \left(\frac{\pi(x+y)}{\pi(x)} - 1\right) q(y) \, d\mu_d(y).$$

Set $l(x) := -\log \pi(x)$, $R(V, x, y) := V(x+y) - V(x) + z V(x)\langle \nabla l(x), y \rangle$ and $R(\pi, x, y) := \pi(x+y)/\pi(x) - 1 + \langle \nabla l(x), y \rangle$. It is proved in Lemma B.4 of [Fort and Moulines (2000)](#) that there exists a constant $c$ such that, for large $|x|$,

$$(3.7) \qquad \sup_{|y| \le c(q)} |R(\pi, x, y)||y|^{-2} \le c|x|^{2(m-1)}.$$

Using a Taylor expansion with integral remainder term of the function $x \mapsto V(x)$, it is easily shown that there exists a constant $c$ such that, for all $z \in (0, z_0)$ and large $|x|$,

$$(3.8) \qquad \sup_{|y| \le c(q)} |R(V, x, y)||y|^{-2} \le c z^2 V(x)|x|^{2(m-1)}.$$

Since $q$ is symmetric, we have

$$PV(x) - V(x)$$
$$= -z V(x) \int_{\mathcal{R}(x)} \langle \nabla l(x), y \rangle^2 q(y) \, d\mu_d(y)$$
$$+ \int R(V, x, y) q(y) \, d\mu_d(y)$$
$$- \int_{\mathcal{R}(x)} R(V, x, y)\langle \nabla l(x), y \rangle q(y) \, d\mu_d(y)$$
$$+ z V(x) \int_{\mathcal{R}(x)} \langle \nabla l(x), y \rangle R(\pi, x, y) q(y) \, d\mu_d(y)$$
$$+ \int_{\mathcal{R}(x)} R(V, x, y) R(\pi, x, y) q(y) \, d\mu_d(y).$$

Thus, for large $|x|$, we deduce from (3.7) and (3.8) that

$$\frac{PV(x) - V(x)}{V(x)}$$
$$\le -z \int_{\mathcal{R}(x)} \langle \nabla l(x), y \rangle^2 q(y) \, d\mu_d(y) + c z^2 |x|^{2(m-1)},$$

for some positive constant $c$ that does not depend on $z$. It is shown in Lemma B.3 of Fort and Moulines (2000) that there exists $\eta > 0$ such that, for large $|x|$,

$$(3.9) \qquad \int_{\mathcal{R}(x)} \langle \nabla l(x), y \rangle^2 q(y) \, d\mu_d(y) > \eta |\nabla l(x)|^2 > \eta d_1^2 |x|^{2(m-1)}.$$

Hence, upon noting that $z d_0 |x|^m \leq \log V(x)$, there exists a constant $\kappa$ which is positive for $z$ small enough, such that, for large $|x|$,

$$PV(x) - V(x) \leq -\kappa [\log V(x)]^{-2(1-m)/m} \, V(x).$$

Since $\pi$ is bounded on compact sets, $\sup_{|x| \leq M} PV(x) + V(x) < \infty$ and the proof is concluded. $\quad \square$

3.3. *Nonlinear autoregressive model.* Consider a process $(\Phi_n, n \geq 0)$ that satisfies the following nonlinear autoregressive equation of order 1:

$$(3.10) \qquad \Phi_{n+1} = g(\Phi_n) + \varepsilon_{n+1},$$

where the sequence $(\varepsilon_n, n \leq 0)$ and the function $g$ satisfy the following assumptions.

ASSUMPTION 3.3. $(\varepsilon_n, n \geq 0)$ is a sequence of i.i.d. zero mean, $d$-dimensional random vectors, independent of $\Phi_0$, that satisfy

$$(3.11) \qquad \mathbb{E}\left[e^{z_0 |\varepsilon_0|^{\gamma_0}}\right] < \infty,$$

for some $z_0 > 0$ and $\gamma_0 \in (0,1]$ and the distribution of $\varepsilon_0$ has a nontrivial absolutely continuous component which is bounded away from zero in a neighborhood of the origin.

ASSUMPTION 3.4. $g : \mathbb{R}^d \to \mathbb{R}^d$ is continuous, and there exist $r, R_0 > 0$ and $\rho \in [0,2)$ such that

$$(3.12) \qquad |g(x)| \leq |x|(1 - r|x|^{-\rho}) \qquad \text{if } |x| \geq R_0.$$

There already exists a wide literature on conditions implying a geometric rate of convergence for nonlinear autoregressive models [see, e.g., Duflo (1997) and Grunwald, Hyndman, Tedesco and Tweedie (2000) and the references therein]. Conditions implying a polynomial rate of convergence have been obtained by Tuominen and Tweedie (1994) and AngoNze (1994) and have been refined by Veretennikov (1997, 2000), AngoNze (2000) and Fort and Moulines (2003). Conditions implying a truly subexponential rate of convergence are considered in Klokov and Veretennikov (2002) [see also Malyshkin (2001) for diffusion processes].

THEOREM 3.3.  *Assume that Assumptions 3.3 and 3.4 hold.*

(i) *If $\rho > \gamma_0$, the drift condition $\mathbf{D}(\phi, V, C)$ holds with $\phi(v) := cv(1 + \log(v))^{1-\rho/(\gamma_0 \wedge (2-\rho))}$, $V(x) := e^{z|x|^{\gamma_0 \wedge (2-\rho)}}$ and $C := \{x \in \mathbb{R}^d, |x| \leq M\}$ for some $z \in (0, z_0)$, $c > 0$ and $M \geq R_0$.*

(ii) *If $\rho \leq \gamma_0$, then the Foster–Lyapunov condition (1.2) holds with $C = \{x \in \mathbb{R}^d, \ |x| \leq M\}$ for some $M \geq R_0$ and $V(x) = e^{z|x|^{\gamma_0}}$ with $z = z_0$ if $\rho < \gamma_0$ and $z \in (0, z_0)$ if $\rho = \gamma_0$.*

COROLLARY 3.4.  *Assume in addition that, for all $x \in \mathbb{R}^d$, $|g(x)| \leq |x|$. Then the chain is $\psi$-irreducible and aperiodic and compact sets of $\mathcal{B}^+(\mathbb{R}^d)$ are petite.*

*If $\rho > \gamma_0$, then there exists $c > 0$ and $z \in (0, z_0)$ such that any probability measure $\lambda$ on $\mathbb{R}^d$ satisfying $\lambda(V) < \infty$ is $(f, r)$-regular with $r(n) = e^{cn^{\{\gamma_0 \wedge (2-\rho)\}/\rho}}$ and $f(x) = e^{z|x|^{\gamma_0 \wedge (2-\rho)}}$.*

PROOF OF THEOREM 3.3.  Throughout the proof, $c$ is a generic constant that can change upon each appearance. Applying the inequality $V(u+w) \leq V(u)V(w)$, we obtain that in all cases, $PV$ is bounded on compact sets of $\mathbb{R}^d$. Thus the proof consists in bounding $PV - V$ outside balls.

(i) We start by examining the case $\rho > \gamma_0$. Set $\beta = \gamma_0 \wedge (2 - \rho)$. We write

$$(3.13) \qquad \frac{PV(x)}{V(x)} - 1 = \frac{PV(x) - V(g(x))}{V(x)} + \frac{V(g(x))}{V(x)} - 1.$$

Using the inequality $(1-u)^{\gamma_0} \leq 1 - \gamma_0 u$ for all $0 \leq u \leq 1$, we have for $|x| \geq R_0$, $|g(x)|^\beta \leq |x|^\beta - \beta r |x|^{\beta-\rho}$, and since $e^x - 1 \leq x + x^2/2$ for all $x \leq 0$,

$$(3.14) \qquad \begin{aligned} \frac{V(g(x))}{V(x)} - 1 &= e^{z|g(x)|^\beta - z|x|^\beta} - 1 \\ &\leq -zr\beta |x|^{\beta-\rho} + \frac{1}{2} z^2 r^2 \beta^2 |x|^{2(\beta-\rho)}. \end{aligned}$$

Let $0 < \eta < 1$. We establish that, for small $z$, large $|x|$ and large $|g(x)|$,

$$(3.15) \qquad PV(x) - V(g(x)) \leq \tfrac{1}{2} z^2 \beta^2 \mathbb{E}[|\varepsilon_0|^2 V(\varepsilon_0)] |x|^{2\beta-2} V(x).$$

Set $R(u, w) = V(u+w) - V(u) - \langle \nabla V(u), w \rangle$. Since $\mathbb{E}[\varepsilon_0] = 0$, this yields

$$(3.16) \qquad \begin{aligned} PV(x) - V(g(x)) &= \mathbb{E}[V(g(x) + \varepsilon_0)] - V(g(x)) \\ &= \mathbb{E}[R(g(x), \varepsilon_0)], \end{aligned}$$

and we have to upper bound the remainder term $\mathbb{E}[R(g(x), \varepsilon_0)]$. If $|w| \leq \eta|u|$, then by using a Taylor expansion with integral remainder term, one has

$$\begin{aligned} |R(u, w)| &\leq \int_0^1 (1-t)|w'\nabla^2 V(u+tw)w| \, dt \\ &\leq \tfrac{1}{2}|w|^2 z\beta \sup_{t \in [0,1]} \{(1 + z\beta|u+tw|^\beta)|u+tw|^{\beta-2} V(u+tw)\}. \end{aligned}$$

Since $y \mapsto |y|^{2\beta-2}e^{z|y|^\beta}$ and $y \mapsto |y|^{\beta-2}e^{z|y|^\beta}$ are ultimately nondecreasing, for large $|u|$ and $|w| \leq \eta|u|$, we have

$$\begin{aligned}
(3.17) \quad |R(u,w)| &\leq \tfrac{1}{2}|w|^2 z\beta\big(1 + z\beta(|u| + |w|)^\beta\big)(|u| + |w|)^{\beta-2}V(u)V(w) \\
&\leq \tfrac{1}{2}z^2\beta^2|w|^2V(w)|u|^{2\beta-2}V(u) + c|w|^2V(w)|u|^{\beta-2}V(u).
\end{aligned}$$

If $|w| \geq \eta|u|$, using again the inequality $V(u+w) \leq V(u)V(w)$,

$$\begin{aligned}
(3.18) \quad |R(u,w)| &\leq V(u+w) + V(u) + |\nabla V(u)||w| \\
&\leq c|w|V(w)|u|^{\beta-1}V(u) \\
&\leq c|w|^2V(w)|u|^{\beta-2}V(u).
\end{aligned}$$

We now apply (3.17) and (3.18) with $u = g(x)$ and $w = \varepsilon_0$. Since $y \mapsto |y|^{2\beta-2}e^{z|y|^\beta}$ and $y \mapsto |y|^{\beta-2}e^{z|y|^\beta}$ are ultimately nondecreasing, for large $|g(x)|$, we have

$$\begin{aligned}
(3.19) \quad &|R(g(x), \varepsilon_0)| \\
&\leq \tfrac{1}{2}z^2\beta^2|\varepsilon_0|^2V(\varepsilon_0)|x|^{2\beta-2}V(x) + c|\varepsilon_0|^2V(\varepsilon_0)|x|^{\beta-2}V(x).
\end{aligned}$$

Equation (3.15) now follows from (3.19). Gathering (3.14) and (3.15), as $\beta \leq 2 - \rho$, we obtain that we can choose $z < z_0$, $M_1$ and $M_2$ such that, for $|x| \geq M_1$ and $|g(x)| \geq M_2$, it holds that

$$PV(x) - V(x) \leq \phi(V(x))$$

with $\phi(v) = -\kappa\beta z\rho^{\rho/\beta}\{1 + \log(v)\}^{1-\rho/\beta}v$ and

$$\kappa = \begin{cases} r, & \text{if } \beta < 2 - \rho, \text{ that is, } \gamma_0 < 2 - \rho, \\ r - 1/2\beta z\mathbb{E}\big[\varepsilon_0^2 e^{z|\varepsilon_0|^\beta}\big], & \text{if } \beta = 2 - \rho, \text{ that is, } \gamma_0 \geq 2 - \rho, \end{cases}$$

and $z$ is chosen small enough such that $\kappa > 0$. To conclude, note that if $|g(x)| \leq M_2$, then $PV(x) \leq V(M_1)\mathbb{E}[V(\varepsilon_0)]$. Choose $M_1$ such that if $|x| \geq M_1$, then $\phi(V(x)) \geq V(M_1)\mathbb{E}[V(\varepsilon_0)]$. Then, defining $C = \{|x| > M_1\}$, we have that, for all $x \notin C$, $PV(x) - V(x) \leq \phi(V(x))$.

(ii) We now consider the case $\rho = \gamma_0$ [observe that $\beta := \gamma_0 \wedge (2 - \rho) = \gamma_0$ and that many results above remain valid]. By (3.13), (3.14), (3.16) and (3.19), we have for large $|x|$ and large $|g(x)|$,

$$\begin{aligned}
\frac{PV(x) - V(x)}{V(x)} &\leq -zr\gamma_0 + \frac{1}{2}z^2r^2\gamma_0^2 \\
&\quad + \frac{1}{2}z^2\gamma_0^2|x|^{2\gamma_0-2}\mathbb{E}[\varepsilon_0^2V(\varepsilon_0)](1 + o(1)).
\end{aligned}$$

For $z$ small enough, the term on the right-hand side is in the interval $(-1, 0)$ and this shows that the Foster–Lyapunov drift condition (1.2) holds with $C$ of the form $\{x, |g(x)| \leq M_1\} \cup \{x, |x| \leq M_2\}$ for large enough $M_1$, $M_2$.

(iii) We finally consider the case $\rho < \gamma_0$. Using the inequality $(1 - u)^{\gamma_0} \leq 1 - \gamma_0 u$ for all $0 \leq u \leq 1$, we have for $|x| \geq R_0$, $|g(x)|^{\gamma_0} \leq |x|^{\gamma_0} - \gamma_0 r |x|^{\gamma_0 - \rho}$. Hence, since $V(u + w) \leq V(u)V(w)$, this yields, for $|x| \geq R_0$,

$$PV(x) = \mathbb{E}\big[ V\big(g(x) + \varepsilon_0\big)\big]$$
$$\leq V(g(x))\mathbb{E}[V(\varepsilon_0)]$$
$$\leq e^{-r\gamma_0 z_0 |x|^{\gamma_0 - \rho}} \mathbb{E}\Big[e^{z_0 |\varepsilon_0|^{\gamma_0}}\Big] V(x).$$

Hence $\lim_{|x| \to \infty} PV(x)/V(x) = 0$, which implies that the Foster–Lyapunov drift condition (1.2) holds with $C := \{|x| \leq M\}$ for large enough $M$.  □

3.4. *Stochastic unit root.* We now consider a process which belongs to the wide family of stochastic unit root models. See, for example, Granger and Sawnson (1997) for many examples. The model we consider is one of the simplest. It has been considered in Gourieroux and Robert (2001) with main focus on its extremal behavior:

$$(3.20) \qquad \mathbf{\Phi}_{n+1} = \mathbf{1}_{\{U_{n+1} \leq g(\mathbf{\Phi}_n)\}} \mathbf{\Phi}_n + \varepsilon_{n+1},$$

where $(\varepsilon_n, n \in \mathbb{N})$ is a sequence of i.i.d. random variables that satisfies (3.11) and $(U_n, n \geq 1)$ is a sequence of i.i.d. random variables, uniformly distributed on $[0,1]$ and independent from the sequence $(\varepsilon_n, n \in \mathbb{N})$. Moreover, we make the following assumption on $g$.

ASSUMPTION 3.5. $g$ is a continuous function with values in $[0,1)$ and there exist $\kappa \in (0,1)$, $c_+(g) > 0$, $c_-(g) < 1$ and $R_0 > 0$ such that

$$(3.21) \qquad \forall x \geq R_0 \qquad 1 - g(x) \geq c_+(g)x^{-\kappa},$$

$$(3.22) \qquad \forall x \leq R_0 \qquad g(x) \leq c_-(g).$$

Let $P$ be the transition kernel of the chain. For all $x \in \mathbb{R}$ and all Borel sets $A$, it can be expressed as

$$(3.23) \qquad P(x, A) = g(x)\mathbb{P}(x + \varepsilon_0 \in A) + (1 - g(x))\mathbb{P}(\varepsilon_0 \in A).$$

Under Assumption 3.5, for all $M > 0$, there exists a constant $\eta(M)$ such that, for all $x \leq M$ and all Borel sets $A$,

$$(3.24) \qquad P(x, A) \geq \eta(M)\mathbb{P}(\varepsilon_0 \in A).$$

This means that every set of the form $(-\infty, M]$ is petite. Define $x_+ = \max(x, 0)$.

THEOREM 3.5. *Under Assumption 3.5 and if $\varepsilon_0$ satisfies (3.11), there exist $z \in (0, z_0]$, $\delta > 0$ and $M \geq R_0$ such that the drift condition $\mathbf{D}(\phi, V, C)$ holds with $V(x) = e^{z x_+^{\beta}}$, $\phi(v) = \delta z^{\tau/\beta} v \{1 \vee \log(v)\}^{-\tau/\beta}$, $C = (-\infty, M]$ and $\beta$ and $\tau$ are given according to the value of $\mathbb{E}[\varepsilon_0]$ by:*

(i) $\beta = \gamma_0 \wedge (1 - \kappa)$ and $\tau = \kappa$, if $\mathbb{E}[\varepsilon_0] > 0$;

(ii) $\beta = \gamma_0 \wedge (1 - \kappa/2)$ $\tau = \kappa$, if $\mathbb{E}[\varepsilon_0] = 0$;

(iii) $\beta = \gamma_0$ and $\tau = (1 - \gamma_0) \wedge \kappa$, if $\mathbb{E}[\varepsilon_0] < 0$.

COROLLARY 3.6.  *Under the same assumptions, the chain is $\psi$-irreducible and (strongly) aperiodic and there exist $c > 0$ and $z > 0$ such that any probability measure $\lambda$ on $\mathbb{R}^d$ satisfying $\lambda(V) < \infty$ is $(f, r)$-regular with:*

(i) $r(n) = e^{cn^{(\gamma_0 \wedge (1-\kappa))/(\gamma_0 \wedge (1-\kappa) + \kappa)}}$ and $f(x) = e^{zx_+^{\gamma_0 \wedge (1-\kappa)}}$, if $\mathbb{E}[\varepsilon_0] > 0$;

(ii) $r(n) = e^{cn^{(\gamma_0 \wedge (1-\kappa/2))/(\gamma_0 \wedge (1-\kappa/2) + \kappa)}}$ and $f(x) = e^{zx_+^{\gamma_0 \wedge (1-\kappa/2)}}$, if $\mathbb{E}[\varepsilon_0] = 0$;

(iii) $r(n) = e^{cn^{\gamma_0/(\kappa \wedge (1-\gamma_0) + \gamma_0)}}$ and $f(x) = e^{zx_+^{\gamma_0}}$, if $\mathbb{E}[\varepsilon_0] < 0$.

PROOF OF THEOREM 3.5.  Let $z < z_0$ and $x > 0$. Using the definition of the transition kernel $P$, we have

$$PV(x) - V(x) = g(x)\mathbb{E}[V(x + \varepsilon_0)] + (1 - g(x))\mathbb{E}[V(\varepsilon_0)] - V(x)$$
$$= g(x)\big(\mathbb{E}[V(x + \varepsilon_0)] - V(x)\big) - (1 - g(x))\big(V(x) - \mathbb{E}[V(\varepsilon_0)]\big)$$
$$\leq \mathbb{E}[V(x + \varepsilon_0)] - V(x) - (1 - g(x))\big(V(x) - \mathbb{E}[V(\varepsilon_0)]\big).$$

Define $R(x, \varepsilon_0) = V(x + \varepsilon_0) - V(x) - \varepsilon_0 \beta z x^{\beta - 1} V(x)$. For any $\eta \in (0, 1)$, we can write

$$\mathbb{E}[V(x + \varepsilon_0)] - V(x) - \beta z \mathbb{E}[\varepsilon_0] x^{\beta - 1} V(x)$$
$$= \mathbb{E}\big[R(x, \varepsilon_0)\mathbf{1}_{\{|\varepsilon_0| \leq \eta x\}}\big] + \mathbb{E}\big[R(x, \varepsilon_0)\mathbf{1}_{\{|\varepsilon_0| > \eta x\}}\big].$$

By the same arguments as in the proof of Theorem 3.3, we have

$$(3.25) \quad \begin{aligned} &\mathbb{E}\big[R(x, \varepsilon_0)\mathbf{1}_{\{|\varepsilon_0| > \eta x\}}\big] \\ &\leq \mathbb{E}\big[V\big((1 + \eta^{-1})|\varepsilon_0|\big) + V(|\varepsilon_0|) + \beta z \eta^{1-\beta}|\varepsilon_0|^\beta V(|\varepsilon_0|)\big]. \end{aligned}$$

Thus this term is bounded provided that $\eta$ and $z$ are chosen such that $(1 + \eta^{-1})^\beta z \leq z_0$. To bound the second term, note that for large enough $x$, the function $x \mapsto x^{2\beta - 2} V(x)$ is increasing. Thus, for $x \geq M$, for some $M$ depending on $\eta$, and $|\varepsilon_0| \leq \eta x$, there exists $t \in (0, 1)$ such that

$$V(x + \varepsilon_0) - V(x) - \beta z \varepsilon_0 x^{\beta - 1} V(x)$$
$$= \tfrac{1}{2}\beta(\beta - 1)z(x + t\varepsilon_0)^{\beta - 2}\varepsilon_0^2 V(x + t\varepsilon_0)$$
$$\quad + \tfrac{1}{2}\big(\beta z(x + t\varepsilon_0)^{\beta - 1}\big)^2 \varepsilon_0^2 V(x + t\varepsilon_0)$$
$$\leq \tfrac{1}{2}\beta^2 z^2 (1 + \eta)^{2\beta - 2} x^{2\beta - 2} \varepsilon_0^2 V(x) V(|\varepsilon_0|)$$
$$\leq \tfrac{1}{2}\beta^2 z^2 x^{2\beta - 2} \varepsilon_0^2 V(x) V(|\varepsilon_0|).$$

For $c < c_+(g)$ and $x$ large enough, say $x \geq M$ for some $M \geq R_0$, we have

$$\big(1 - g(x)\big)\big(V(x) - \mathbb{E}[V(\varepsilon_0)]\big) \geq cx^{-\kappa}V(x).$$

Hence, taking (3.25) into account, there exists a positive real number $M$ such that, if $x \geq M$, then

$$PV(x) - V(x)$$
$$\leq \big(z\beta x^{\beta-1}\mathbb{E}[\varepsilon_0] + \tfrac{1}{2}\beta^2 z^2 x^{2\beta-2}\mathbb{E}[\varepsilon_0^2 V(|\varepsilon_0|)] - cx^{-\kappa}\big)V(x).$$

If $\mathbb{E}[\varepsilon_0] > 0$, set $\beta = \gamma_0 \wedge (1 - \kappa)$. Then, for large enough $x$, we obtain

$$PV(x) - V(x) \leq -\delta x^{-\kappa}V(x)$$
$$= -\delta z^{\kappa/\beta}V(x)\{\log(V(x))\}^{-\kappa/\beta},$$

with $\delta = c < c_+(g)$ if $\gamma_0 < 1 - \kappa$ or $\delta = c - \beta z\mathbb{E}[\varepsilon_0]$, $c < c_+(g)$ and $z$ such that $\delta > 0$ if $\gamma_0 \geq 1 - \kappa$. If $\mathbb{E}[\varepsilon_0] < 0$, set $\beta = \gamma_0$ and $\tau = (1 - \gamma_0) \wedge \kappa$. Then, for $x$ large enough,

$$PV(x) - V(x) \leq -\delta x^{-\tau}V(x)$$
$$= -\delta z^{\tau/\gamma_0}V(x)\{\log(V(x))\}^{-\tau/\gamma_0},$$

with $\delta = c < c_+(g)$ if $\gamma_0 < 1 - \kappa$ and $\delta = c - z\beta\mathbb{E}[\varepsilon_0]$, $c < c_+(g)$ and $z$ such that $\delta > 0$ if $\gamma_0 \geq 1 - \kappa$. If $\mathbb{E}[\varepsilon_0] = 0$, then $\beta$ must satisfy $2\beta - 2 \leq -\kappa$; thus we set $\beta = (1 - \kappa/2) \wedge \gamma_0$, and we obtain

$$PV(x) - V(x) \leq -\delta x^{-\kappa}V(x)$$
$$= -\delta z^{\kappa/\beta}V(x)\{\log(V(x))\}^{-\kappa/\beta},$$

with $\delta = c < c_+(g)$ if $1 - \kappa/2 > \gamma_0$ and $\delta = c - \tfrac{1}{2}\beta^2 z^2 \mathbb{E}[\varepsilon_0^2 V(|\varepsilon_0|)]$, with $c < c_+(g)$ and $z$ such that $\delta > 0$ if $1 - \kappa/2 \leq \gamma_0$. $\square$

## REFERENCES

ANGONZE, P. (1994). Critres d'ergodicité de modèles markoviens. Estimation non-paramétrique sous des hypothèses de dépendance. Ph.D. dissertation, Univ. Paris 9, Dauphine.

ANGONZE, P. (2000). Geometric and subgeometric rates for Markovian processes: A robust approach. Technical report, Univ. de Lille III.

DUFLO, M. (1997). *Random Iterative Systems.* Springer, Berlin. MR1485774

FORT, G. and MOULINES, E. (2000). V-subgeometric ergodicity for a Hastings–Metropolis algorithm. *Statist. Probab. Lett.* **49** 401–410. MR1796485

FORT, G. and MOULINES, E. (2003). Polynomial ergodicity of Markov transition kernels. *Stochastic Process. Appl.* **103** 57–99. MR1947960

GOURIEROUX, C. and ROBERT, C. (2001). Tails and extremal behaviour of stochastic unit root models. Technical report, Centre de Recherche en Economie et Statistique du Travail.

GRANGER, C. and SAWNSON, N. (1997). An introduction to stochastic unit-root processes. *J. Econometrics* **80** 35–62. MR1467611

Grunwald, G., Hyndman, R., Tedesco, L. and Tweedie, R. (2000). Non-Gaussian conditional linear AR(1) models. *Aust. N. Z. J. Stat.* **42** 479–495. MR1802969

Jarner, S. and Hansen, E. (2000). Geometric ergodicity of Metropolis algorithms. *Stochastic Process. Appl.* **85** 341–361. MR1731030

Jarner, S. and Roberts, G. (2002). Polynomial convergence rates of Markov chains. *Ann. Appl. Probab.* **12** 224–247. MR1890063

Klokov, S. and Veretennikov, A. (2002). Sub-exponential mixing rate for a class of Markov processes. Technical Report 1, School of Mathematics, Univ. Leeds.

Malyshkin, M. (2001). Subexponential estimates of the rate of convergence to the invariant measure for stochastic differential equations. *Theory Probab. Appl.* **45** 466–479. MR1967786

Mengersen, K. and Tweedie, R. (1996). Rates of convergence of the Hastings and Metropolis algorithms. *Ann. Statist.* **24** 101–121. MR1389882

Meyn, S. and Tweedie, R. (1993). *Markov Chains and Stochastic Stability.* Springer, London. MR1287609

Nummelin, E. and Tuominen, P. (1983). The rate of convergence in Orey's theorem for Harris recurrent Markov chains with applications to renewal theory. *Stochastic Process. Appl.* **15** 295–311. MR711187

Roberts, G. and Tweedie, R. (1996). Geometric convergence and central limit theorem for multidimensional Hastings and Metropolis algorithms. *Biometrika* **83** 95–110. MR1399158

Tanikawa, A. (2001). Markov chains satisfying simple drift conditions for subgeometric ergodicity. *Stoch. Model.* **17** 109–120. MR1853436

Tuominen, P. and Tweedie, R. (1994). Subgeometric rates of convergence of $f$-ergodic Markov chains. *Adv. in Appl. Probab.* **26** 775–798. MR1285459

Veretennikov, A. (1997). On polynomial mixing bounds for stochastic differential equations. *Stochastic Process. Appl.* **70** 115–127. MR1472961

Veretennikov, A. (2000). On polynomial mixing and convergence rate for stochastic differential and difference equations. *Theory Probab. Appl.* **44** 361–374. MR1751475

R. Douc
CMAP
Ecole Polytechnique
Palaiseau
France

E. Moulines
Département TSI
Ecole Nationale Supérieure
    des Télécommunications
46 rue Barrault
75013 Paris
France

G. Fort
LMC-IMAG
51, rue des Mathématiques
BP53, 38041 Grenoble Cedex 9
France

P. Soulier
Département de Mathématiques
Université d'Evry
91025 Evry Cedex
France
e-mail: phlippe.soulier@maths.univ-evry.fr